\documentclass[12pt]{article}

\usepackage{graphicx}

\usepackage{amsmath}
\usepackage{amssymb}
\usepackage{color}
\usepackage{enumerate}

\def\Chi{\raise .3ex \hbox{\large $\chi$}}

\def\cO{{\cal O}}

\def\R{\mathbb{R}} 
\def\Z{\mathbb{Z}} 
\def\N{\mathbb{N}}

\newcommand{\cF}{\mathcal{F}}

\def\dist{{\rm dist}}

\def\cH{{\cal H}}

\def\span{{\rm span}}

\newcommand{\rr}{\rangle}
\renewcommand{\lll}{\langle}

\newtheorem{lemma}{Lemma}[section]

\newtheorem{cor}[lemma]{Corollary}
\newtheorem{theorem}[lemma]{Theorem}
\newtheorem{rem}[lemma]{Remark}

\numberwithin{equation}{section} 
\newcommand{\beqn}{\begin{equation}}
\newcommand{\eeqn}{\end{equation}}
\newcommand{\be}{\begin{equation}}
\newcommand{\ee}{\end{equation}}

\newcommand\eref[1]{(\ref{#1})}

\def\dist{\mathrm{dist}}

\def\rnew{\color{red}}

\title{Greedy Algorithms for Reduced Bases in Banach Spaces\thanks{%
 This research was supported by
the Office of Naval Research Contracts
   ONR-N00014-08-1-1113, 
   ONR N00014-09-1-0107, and
   ONR N00014-11-1-0712;
   the AFOSR Contract FA95500910500;
the NSF Grants
   DMS-0810869, and
   DMS 0915231; and 
the  EU Project POWIEW.  This publication is based in part on work supported by Award No. KUS-C1-016-04  made by King Abdullah University of Science and Technology (KAUST)}
\author{
 Ronald DeVore, Guergana Petrova, and Przemyslaw
Wojtaszczyk}
} 

\begin{document}

\maketitle

\begin{abstract}   Given a Banach space $X$ and one of its compact sets $\cF$, we consider the problem of finding a  good $n$ dimensional space $X_n\subset X$  which can be used to approximate the elements of $\cF$.    The best possible error we can achieve for such an approximation is given by the  Kolmogorov width $d_n(\cF)_X$.   However, finding the space which  gives this performance is typically numerically intractable.   Recently, a new greedy strategy for obtaining good spaces was given  in the context of  the reduced basis method for solving a parametric family of PDEs.   The performance of this greedy algorithm was initially analyzed in  \cite{BMPPT} in the case $X=\cH$ is a Hilbert space. The results of \cite{BMPPT} were significantly  improved on in \cite{BCDDPW}.   The purpose of the present paper is to give a new analysis of  the performance of such greedy algorithms. Our analysis not only gives improved results for the Hilbert space case but  can also  be applied to the same greedy procedure in general Banach spaces.\end{abstract}

\section{Introduction} \label{intro}

 Let $X$ be a Banach space with norm $\|\cdot\|:=\|\cdot\|_X$, and let $\cF$ be one of its   compact subsets.  For notational convenience only,
we shall assume that the elements $f$ of $\cF$ satisfy $\|f\|_X\le 1$. We consider the following greedy algorithm for generating approximation spaces for $\cF$. We first choose a function $f_0$ such that
\be
\label{first}
{\displaystyle \|f_0\| =\max_{f\in\cF}\|f\|.}
 \ee
  Assuming $\{f_0,\dots,f_{n-1}\}$ and $V_n:=\span\{f_0,\dots,f_{n-1}\}$ have been selected, we then take $f_n\in\cF$ such that   
\be
\label{second}
 \dist( f_n,V_n)_X\|=\displaystyle{\max_{f\in \cF} \dist(f,V_n)_X},
\ee
  and define
  \be
  \label{second1}
   \sigma_n:=\sigma_ n(\cF)_X:=\dist( f_n,V_n)_X:=\sup_{f\in\cF}\inf_{g\in V_n}\|f-g\|.
  \ee
  This greedy algorithm was introduced, for  the case $X$ is  a Hilbert space, in the reduced basis method \cite {MPT, MPT1} for solving a family of PDEs.  Certain variants of this algorithm,  known as {\it weak greedy algorithms},   described below,  are now numerically implemented with great success  in   the reduced basis method.   Our interest in this paper will be in the approximation properties of this algorithm and its weak variant.

   We are interested in how well the space $V_n$  approximates  the elements of $\cF$ and for this purpose we compare its performance
   with the best possible performance which is given by the Kolmogorov width $d_n(\cF)_ X$ of $\cF$ defined by

\beqn
\label{kolw}
d_n:=d_n(\cF)_X:=\inf_{Y} ~ \sup_{f\in \cF} ~~ \dist(f,Y)_X, 
\eeqn
  where the infimum is taken over all $n$ dimensional  subspaces $Y$ of $X$. 
We refer  the reader to \cite{LGM} for a general discussion of Kolmogorov widths. We also define
\beqn
\nonumber
d_0:=d_0(\cF)_X:=\max_{f\in\cF}\|f\|=\sigma_0(\cF)_X,
\eeqn
which corresponds to approximating by zero dimensional spaces.

Of course, if $(\sigma_n)_{ n\geq 0}$  decays  at a rate comparable to $(d_n)_{n\geq 0}$, this would mean that  the greedy selection provides essentially the best possible accuracy attainable by $n$-dimensional subspaces. 
Various  comparisons have been given between $\sigma_n$ and $d_n$.    A first result in this direction, in the case that $X$ is a Hilbert space $\cH$,   was   given in \cite{BMPPT}  where it was proved that%
\beqn
\label{BMPPT}
\sigma_n(\cF)_\cH\le Cn2^nd_n(\cF)_\cH, 
\eeqn
with $C$ an absolute constant.
While this is an interesting comparison, it is only useful if $d_n(\cF)_\cH$ decays to zero faster than $n^{-1}2^{-n}$.  
 Various improvements on \eref{BMPPT} were given in \cite{BCDDPW}, again in the Hilbert space setting.   We mention two of these.  It was shown that if $d_n(\cF)_\cH\le Cn^{-\alpha}$, $n=1,2,\dots$, then
 \be
 \label{poly1}
 \sigma_n(\cF)_\cH\le C_\alpha'n^{-\alpha}.
 \ee
 This shows that in the scale of polynomial decay the greedy algorithm performs with the same rates as $n$-widths.
 A related result was proved for sub-exponential decay.    If for some $0<\alpha\le 1$, we have  
 $d_n(\cF)_\cH\le C e^{- c n^\alpha}$,  $n=1,2,\dots$, then
 \be
 \label{poly2}
 \sigma_n(\cF)_\cH\le C_\alpha 'e^{-c_{\alpha}'n^{\beta}},\quad {\beta}=\frac{\alpha}{\alpha +1}, \quad n=1,2,\dots.
 \ee

 In numerical implementations, the greedy algorithm  is too demanding since at each iteration it requires finding an element from $\cF$ which is at furthest distance from $V_ n$.   To circumvent this difficulty, one modifies the algorithm as follows.  We fix a constant $0<\gamma\le 1$.
  At the first step of the algorithm, 
one  chooses a function $f_0\in\cF$ such that
\beqn
\nonumber
{ \|f_0\|\ge \gamma \sigma_0(\cF)_X.}
\eeqn
 At the general step, if $f_0,\dots,f_{n-1}$ have been chosen, 
 $  V_{n}:=\span\{f_0,\dots,f_{n-1}\}$,  and
\beqn
\nonumber
\sigma_n(f)_X:=\dist(f,V_n)_X,
\eeqn
 is the best approximation error to $f$ from $V_n$  we now
 choose $f_{n}\in \cF$ such that
\beqn
\label{gae1}
{\displaystyle \sigma_n(f_{n})_X\ge \gamma \max_{f\in\cF} \sigma_n(f)_X},
\eeqn
 to be the next element in the greedy selection. Note that if $\gamma=1$, then the weak greedy algorithm 
reduces to the greedy algorithm that we have introduced above.

Notice that similar to the greedy algorithm, $(\sigma_n(\cF)_X)_{n\ge 0}$ is also  monotone decreasing.
Of course, neither the greedy algorithm or the weak greedy algorithm   give a unique sequence $(f_n)_{n\geq 0}$,
nor is the sequence $(\sigma_n(\cF)_X)_{n\ge 0}$ unique.  In all that follows, the notation reflects any sequences which can arise in the implementation of the weak greedy selection for the fixed value of $\gamma$.

  In the present paper, we shall first  
  prove a lemma that we use in our  new analysis of the weak greedy algorithm.
  This new analysis gives a significant improvement of the previous results.   We mention two of these:  
  
  The first, given in  Corollary  \ref{C1},  is that  
  \be
  \label{direct}
  \sigma_{2n}(\cF)_\cH\le \sqrt{2}\gamma^{-1}\sqrt{d_n(\cF)_\cH},\quad n=1,2,\dots. 
  \ee
   This is the first direct comparison between $(\sigma_n(\cF)_\cH)_{ n\geq 0}$ and $(d_n(\cF)_\cH)_{ n\geq 0}$ 
    for the special case $X$ is a Hilbert space $\cH$, which guarantees a specific rate of decay for $(\sigma_n(\cF)_\cH)_{ n\geq 0}$ without any assumption of a decay rate for $(d_n(\cF)_\cH)_{ n\geq 0}$.    Notice that, in particular, this allows one to improve the sub-exponential results mentioned earlier (see Corollary \ref{C1}).

 The second part of our paper analyzes the performance of the greedy algorithm in  a general Banach space.   We prove estimates for the decay of $(\sigma_n(\cF)_X)_{n\geq 0}$ similar to those in the Hilbert space  case, except that there is a loss of the order $O(\sqrt{n})$.   We give examples which show that this loss in essence cannot be removed.  However, our results for a general Banach space are still not definitive.   For example,
 we have no result of the form \eref{direct} because of the $\sqrt{n}$ factor that appears in our Banach space results.

\section{Main lemma}
\label{sect21}

In this section, we shall prove a lemma for matrices that we employ in our analysis of weak greedy algorithms in both Hilbert and 
Banach spaces.
 \begin{lemma}
 \label{L1}
 Let $G=(g_{i,j})$ be a $K\times K$  lower triangular matrix with rows ${\bf g}_1, \ldots,{\bf g}_K$, $W$ be any $m$ dimensional 
 subspace of $\R^K$,  and $P$ be the orthogonal projection of $\R^K$ onto $W$. Then
 \be
 \label{L11}
  \prod_{i=1}^Kg_{i,i}^2\le  \left\{\frac{1}{m}\sum_{i=1}^K \|P{\bf g}_i\|_{\ell_2}^2\right\}^m\left\{ \frac{1}{K-m}\sum_{i=1}^K \|{\bf g}_i-P{\bf g}_i\|_{\ell_2}^2\right\}^{K-m},
 \ee
 where $\|\cdot\|_{\ell_2}$ is the euclidean norm of a vector in $\R^K$.
 \end{lemma}
 
 {\bf Proof:}  We  choose  an orthonormal basis ${\bf \varphi}_1,\dots,{\bf \varphi}_m$ for the space $W$ and complete it into an orthonormal  basis 
 ${\bf \varphi}_1,\dots,{\bf \varphi}_K$ for $\R^K$.  If we denote by $\Phi$ the $K\times K$ orthogonal matrix whose $j$-th column is ${\bf \varphi}_j$, then   
 the matrix $C:=G\Phi$ has entries  $c_{i,j}=\langle {\bf g}_i,{\bf \varphi}_j\rangle$.   
  We denote by ${\bf c}_j$, the $j$-th column of $C$.   It follows from the arithmetic geometric mean inequality 
  for the numbers $\{\|{\bf c}_j\|_{\ell_2}^2\}_{j=1}^m$ that
 \be
 \label{fc}
\prod_{j=1}^m \|{\bf c}_j\|_{\ell_2}^2\le \left\{\frac{1}{m}\sum_{j=1}^m \|{\bf c}_j\|^2_{\ell_2}\right\}^m=
\left\{\frac{1}{m}\sum_{j=1}^m \sum_{i=1}^K\langle {\bf g}_i,\varphi_j\rangle^2\right\}^m=
\left \{\frac{1}{m}\sum_{i=1} ^K\|P{\bf g}_i\|_{\ell_2}^2\right\}^m.
 \ee
 Similarly, 
 
  \begin{eqnarray}
 \label{lc}
 \prod_{j=m+1}^K \|{\bf c}_j\|_{\ell_2}^2&\le& \left\{\frac{1}{K-m} \sum_{j=m+1}^K\|{\bf c}_j\|_{\ell_2}^2\right\}^{K-m}
 = \left\{\frac{1}{K-m} \sum_{i=1}^K\|{\bf g}_i-P{\bf g}_i\|_{\ell_2}^2\right\}^{K-m}.
 \end{eqnarray}
Now,  Hadamard's inequality for the matrix $C$ and relations \eref{fc} and \eref{lc}  result in
\begin{eqnarray}
\label{T13}
(\det C)^2&\le& \prod_{j=1}^K\|{\bf c}_j\|^2_{\ell_2}
\le\left \{\frac{1}{m}\sum_{i=1}^K \|P{\bf g}_i\|_{\ell_2}^2\right\}^m\left\{ \frac{1}{K-m}\sum_{i=1}^K\|{\bf g}_i-P{\bf g}_i\|_{\ell_2}^2\right\}^{K-m}.
\end{eqnarray}%
 The latter inequality and  the fact that $\displaystyle {\det G=\prod_{i=1}^Kg_{i,i}}$ and $|\det C| =|\det G|$ gives \eref{L11}.
$\hfill\Box$

\section{A new analysis for the weak greedy algorithm in a Hilbert space}
\label{sect2}

The purpose of this section is to obtain new results for the performance of the weak greedy algorithm in a Hilbert space that considerably improve
on the analysis in \cite{BMPPT} and \cite{BCDDPW}.  This will be accomplished by making a finer  comparison between $(\sigma_n(\cF)_\cH)_{n\geq 0}$ and $(d_n(\cF)_\cH)_{ n\geq 0}$ than those given in \cite{BCDDPW}.
We assume throughout  this section that $X=\cH$ is a Hilbert space and follow the notation from \cite{BCDDPW}.

 Note that in general, the weak greedy algorithm does not terminate and we obtain an infinite sequence $f_0, f_1,f_2,\dots$.
In order to have a consistent notation in what follows, we shall define $f_m:=0$, $m>N$, if the algorithm terminates at $N$, i.e. if $\sigma_N(\cF)_\cH=0$.
By $(f_n^*)_{n\geq 0}$ we denote the orthonormal system
obtained from  $(f_n)_{n\geq 0}$ by Gram-Schmidt orthogonalization.  It follows that  the orthogonal projector $P_n$ from $\cH$ onto $V_n$ is given by
$$
P_n f=\sum_{i=0}^{n-1}\lll f, f_i^*\rr f_i^*,
$$
and, in particular,
\beqn
\nonumber
f_n = P_{ n+1} f_n =\sum_{j=0}^n a_{n,j} f^*_j,\quad a_{n,j}=\lll f_n,f^*_j\rr, \,\,j\leq n.
\eeqn
There is no loss of generality in assuming that the infinite dimensional Hilbert space $\cH$ is $\ell_2(\N\cup \{0\})$
and that $f_j^*=e_j$, where $e_j$ is the vector with a one in the coordinate indexed by $j$ and  is zero in all other coordinates, i.e. $(e_j)_i =\delta_{j,i}$.  We adhere to this assumption throughout this section of the  paper.

We consider the lower triangular matrix
$$
A:= (a_{i,j})_{i,j=0}^\infty,\quad a_{i,j}:=0,\, j>i.
$$
This matrix incorporates all the   information about the  weak greedy algorithm on $\cF$.
The following two properties characterize any lower triangular matrix $A$ generated by such a greedy algorithm. With the notation $\sigma_n:= \sigma_n(\cF)_\cH$, we have:
\vskip .1in
\noindent
 {\bf P1:}  The diagonal elements of $A$  satisfy
$\gamma \sigma_n \leq |a_{n,n}|\leq \sigma_n$.
\vskip .1in
\noindent
  {\bf P2:} For every $m\ge n$, one has
$\sum_{j=n}^m a_{m,j}^2\leq \sigma_n^2$.
\vskip .1in

\noindent
Indeed, {\bf P1} follows from
$$
a_{n,n}^2 = \|f_n\|^2 -\|P_{n}f_{n}\|^2 =\|f_n-P_{n}f_n\|^2,
$$
combined with the weak greedy selection property \eref{gae1}.
To see {\bf P2}, we note that for  $m\ge n$,
$$
\sum_{j=n}^m a_{m,j}^2=\|f_m-P_{n}f_m\|^2\leq \max_{f\in\cF} \|f-P_{n}f\|^2 =\sigma_n^2.
$$

\begin{rem}
\label{remP1P2}
 If $A$ is any matrix satisfying {\bf P1} and {\bf P2} with $ (\sigma_n)_{n\geq 0}$ a  decreasing sequence that converges to $0$, then
 the rows of $A$ form a compact subset of $\ell_2(\N\cup \{0\})$.  If $\cF$ is the set consisting of these rows, 
 then one of the possible realizations of the weak greedy algorithm  with constant $\gamma$
 will choose the rows in  that order and $A$ will be the resulting  matrix.
 \end{rem}
 
 The matrix representation $A$ of the weak greedy algorithm was the basis of the analysis given in \cite{BCDDPW} and will also be critical in the proof of the next theorem.
 
 \begin{theorem}
 \label{T10}
 For the weak greedy algorithm with constant $\gamma$ in a Hilbert space  $\cH$ and for any compact set $\cF$, we have the following inequalities between $\sigma_n:= \sigma_n(\cF)_\cH$ and $d_n:=d_n(\cF)_\cH$, for any $N\ge 0$, $K\ge 1$, and $1\le m< K$, 
 \be
 \label{T1}
 \prod_{i=1}^K\sigma^2_{N+i} \le \gamma^{-2K}\left\{\frac{K}{m}\right\}^m\left\{\frac{K}{K-m}\right\}^{K-m}
 \sigma_{N+1}^{2m}d_m^{2K-2m}
 \ee

 \end{theorem}
 
 {\bf Proof:}  
  We consider  the $K\times K$ matrix $G=(g_{i,j})$ which is  formed by the rows and columns of $A$ with indices from $\{N+1,\dots,N+K\}$.  
 Each row ${\bf g}_i$ is the restriction of $f_{N+i}$ to the coordinates $N+1,\dots,N+K$.  
Let $\cH_m$ be the $m$-dimensional Kolmogorov subspace of $\cH$  for which $\dist(\cF,\cH_m) = d_m$.  Then, 
$ \dist(f_{N+i},\cH_m) \leq d_m$, $i=1,\ldots K$. Let $\widetilde W$ be the linear space which is the restriction of $\cH_m$ to the coordinates $N+1,\dots,N+K$. In general, $\dim (\widetilde W)\leq m$. Let $W$ be an  $m$ dimensional space, $W\subset \span\{e_{N+1},\dots,e_{N+K}\}$,  such that $\widetilde W\subset W$ and $P$ and $\widetilde P$ are the projections in $\R^K$ onto $W$ and $\widetilde W$, respectively.
 Clearly, 
\be
 \label{also}
 \|P{\bf g}_i\|_{\ell_2}\le \|{\bf g}_i\|_{\ell_2}\le \sigma_{N+1},\quad i=1,\ldots,K,
 \ee 
 where we have used Property  {\bf P2} in the last inequality.
 Note that
 \be
 \label{alsoo}
\|{\bf g}_i-P{\bf g}_i\|_{\ell_2}\leq \|{\bf g}_i-\widetilde P{\bf g}_i\|_{\ell_2}=\dist({\bf g}_i,\widetilde W)\le 
 \dist(f_{N+i},\cH_m) \leq d_m, \quad i=1,\dots , K.
 \ee
 It follows from   Property {\bf P1} that
   \be
 \label{Tha}
 \prod_{i=1}^K|a_{N+i,N+i}|\ge \gamma^K \prod_{i=1}^K\sigma_{N+i}.
 \ee
We now  apply Lemma \ref{L1} for this $G$ and $W$, and use estimates \eref{also}, \eref{alsoo},  and \eref{Tha} to derive \eref{T1}. The proof is 
completed.
$\hfill\Box$

 We next record some special cases of Theorem \ref{T10}.
 \begin{cor}
 \label{C1}
 For the weak greedy algorithm  with constant $\gamma$ in a Hilbert space $\cH$, we have the following:
 
 {\rm (i)} For any compact set $\cF$ and $n\ge 1$, we have
 \be
 \label{C11}
 \sigma_n(\cF)\le  \sqrt{2}\gamma^{-1}\min_{1\le m< n} d_m^{\frac{n-m}{n}}(\cF).
\ee 
In particular $\sigma_{2n}(\cF)\le \sqrt{2}\gamma^{-1}\sqrt{d_n(\cF)}$, $n=1,2\dots$.

{\rm (ii)} If $d_n(\cF)\le C_0n^{-\alpha}$, $n=1,2,\dots$, then   $\sigma_n(\cF)\le C_1n^{-\alpha}$, $n=1,2\dots$, with $C_1:=2^{ 5\alpha+1}\gamma^{-2}C_0$. 

{\rm (iii)} If $d_n(\cF)\le C_0e^{-c_0n^{\alpha}}$, $n=1,2,\dots$, then    $\sigma_n(\cF)\le \sqrt{2C_0}\gamma^{-1}e^{-c_1n^{\alpha}}$, $n=1,2\dots$, where $c_1=2^{-1-2\alpha}c_0$,

 \end{cor}
 
 {\bf Proof:} (i) We take $N=0$, $K=n$ and any $1\le m< n$ in Theorem \ref{T10}, use the monotonicity of 
 $(\sigma_n)_{n\geq 0}$ and the fact that $\sigma_0\le 1$ to  obtain
 \be
 \label{C14}
 \sigma_n^{2n}\le \prod_{j=1}^n \sigma_j^2\le  \gamma^{-2n} \left\{\frac{n}{m}\right \}^m\left\{\frac{n}{n-m}\right\}^{n-m}d_m^{2n-2m}.
 \ee
 Since $x^{-x}(1-x)^{x-1}\le 2$ for $0< x< 1$, we derive \eref{C11}.
 
 (ii)   It follows from the monotonicity of $(\sigma_n)_{ n\geq 0}$ and  \eref{T1} for $N=K=n$ and any $1\leq m<n$
 that
 $$
  \sigma_{2n}^{2n}\le  \prod_{j=n+1}^{2n}\sigma_j^2
\le \gamma^{-2n}\left\{\frac{n}{m}\right\}^m\left\{\frac{n}{n-m} \right\}^{n-m}\sigma_{n}^{2m}d_m^{2n-2m}.
 $$
In the case $n=2s$ and  $m=s$ we have
 \be
 \label{Cnew}
 \sigma_{4s}\leq \sqrt{2}\gamma^{-1}\sqrt{\sigma_{2s}d_s}.
\ee 
Now we prove our claim  by contradiction.  Suppose it is not true and $M$ is  the first value   where $\sigma_M(\cF)>C_1M^{-\alpha}$.    Let us first assume   $M=4s$.   From  \eref{Cnew}, we have 
 \begin{eqnarray}
\label{pp}
\sigma_{4s}\leq  \sqrt{2}\gamma^{-1} \sqrt{C_1(2s)^{-\alpha}}\sqrt{C_0s^{-\alpha}}=\sqrt{2^{1-\alpha}C_0C_1}\gamma^{-1}s^{-\alpha},
  \end{eqnarray}
   where we have used the  fact that $\sigma_{2s}\le C_1(2s)^{-\alpha}$ and $d_s\leq C_0s^{-\alpha}$.  It follows that
 $$
 C_1(4s)^{-\alpha}<\sigma_{4s}\le \sqrt{2^{1-\alpha}C_0C_1}\gamma^{-1}s^{-\alpha},
 $$
 and therefore
 
 \be
  \nonumber
 C_1 < 2^{3\alpha+1}\gamma^{-2}C_0<2^{5\alpha+1}\gamma^{-2}C_0,
 \ee
 which is the desired contradiction.  
 If $M=4s+q$, $q\in\{1,2,3\}$, then it follows from \eref{pp} and the monotonicity of $(\sigma_n)_{ n\geq 0}$  that
 $$
C_12^{-3\alpha} s^{-\alpha}=C_12^{-\alpha}(4s)^{-\alpha}<C_1(4s+q)^{-\alpha}< \sigma_{4s+q}\leq  \sigma_{4s}\le \sqrt{2^{1-\alpha}C_0C_1}\gamma^{-1}s^{-\alpha}.
 $$
 From this, we obtain
 $$
 C_1<2^{5\alpha+1}\gamma^{-2}C_0,
 $$
  which is the desired contradiction in this case.   This completes the proof of (ii).

 (iii) From (i), we have
 \be
 \label{C17}
 \sigma_{2n+1} \le \sigma_{2n} \le \sqrt{2}\gamma^{-1}\sqrt{d_n}\le \sqrt{2C_0}\gamma^{-1}e^{-\frac{c_0}{2}n^\alpha}=\sqrt{2C_0}\gamma^{-1}
 e^{-c_02^{-1-\alpha}(2n)^\alpha},
 \ee
 from which (iii) easily follows.   \hfill $\Box$

 \begin{rem}
 Note that one can obtain a better constant $c_1$ in {\rm  (iii)} if the minimum in {\rm\eref{C11}}  is computed.  Namely, this gives
  \be
 \nonumber
  \sigma_{2n} \le \sqrt{2}\gamma^{-1}C_0\min_{1\le m< n} e^{-c_0m^\alpha\frac{(n-m)}{n}}=
 \sqrt{2}\gamma^{-1}C_0 e^{\displaystyle{-c_0n^\alpha \{\max_{1\le m< n}\left(\frac{m}{n}\right)^\alpha\left(1-\frac{m}{n}\right)\}}}.
 \ee
 Then,  using the fact that  $ x^\alpha(1-x)$, $0<x<1$, has a maximum at $\frac{\alpha}{\alpha+1}$  results in a better constant.
\end{rem}

\section{Bounds for the greedy algorithm in Banach spaces}
 
We will now derive  bounds for the performance of the  weak  greedy algorithm in a general Banach space $X$.  
In this section, we will use the abbreviation $\sigma_n:= \sigma_n(\cF)_X$ and $d_n:=d_n(\cF)_X$.
As in the Hilbert space case, we  associate with the greedy procedure  a lower triangular matrix $A=(a_{i,j})_{i,j=0}^\infty$  in  the following way.
For each $j=0,1,\dots$, we let $\lambda_j\in X^*$ be  the linear functional of norm one  that  satisfies
 \be
 \label{functionals}
  {\rm (i)}\  \lambda_j(V_j)=0, \quad  {\rm (ii)}\ \lambda_j(f_j)=\dist(f_j,V_j)_X. 
\ee
 The existence of such a functional is a simple consequence of the Hahn-Banach theorem (see  \cite[Chapt. IV, Cor.14.13]{W}).   We let $A$ be the matrix
  with  entries 
 $$
 a_{i,j}=\lambda_j(f_i).
 $$
From (ii) of \eref{functionals, we see that $A$ is lower triangular.} Its  diagonal elements $a_{j,j}$ satisfy the inequality
\be
\label{diag}
\gamma \sigma_j\leq  a_{j,j}=\dist(f_j,V_j)_X =\sigma_j,
 \ee
 because of the weak greedy selection property \eref{gae1}.
   Also, each entry $a_{i,j}$ satisfies
 $$
  |a_{i,j}|=|\lambda_j(f_i)| = |\lambda_j(f_i-g)|\leq  \|\lambda_j\|_{X^*}\|f_i-g\|=\|f_i-g\|, \quad j<i,
 $$
 for every $g\in V_j$, since $\lambda_j(V_j)=0$. Therefore we have
 \be
 \label{entries}
 |a_{i,j}|\leq   \dist(f_i,V_{j})\le \sigma_{j}, \quad j<i.
 \ee

\begin{theorem}
 \label{T20}
 For the weak greedy algorithm with constant $\gamma$ in a Banach space $X$ and for any compact set $\cF $ contained in the unit ball of $X$, we have the following inequalities between $\sigma_n:= \sigma_n(\cF)_X$ and $d_n:=d_n(\cF)_X$:  for 
 any $N\ge 0$, $K\ge 1$, and $1\le m < K$, 
 \be
 \label{T3}
 \prod_{i=1}^K\sigma_{N+i}^2\le 2^K K^{K-m}\gamma^{-2K} \left\{ \sum_{i=1}^K \sigma_{N+i}^2\right\}^m d^{2K-2m}_m.
  \ee
 \end{theorem}
 
 {\bf Proof:} 
  As in  the proof of Theorem \ref{T10}, we consider  the   $K\times K$ matrix $G$ which is  formed by the rows and columns of 
 $A$ with indices from $\{N+1,\dots,N+K\}$. 
 Let $X_m$ be  the Kolmogorov subspace of $X$ for which $\dist(\cF,X_m)= d_m$.  For each $i$, there is an element
  $h_i\in X_m$ such that 
  $$
  \|f_i-h_i\|=\dist(f_i,X_m)_X\le d_m,
  $$  
   and therefore
  \be
  \label{approx1}
  |\lambda_j(f_i)-\lambda_j(h_i)|=|\lambda_j(f_i-h_i)|\le \|\lambda_j\|_{X^*}\|f_i-h_i\|\leq d_m.
  \ee

  We now consider the vectors
  $(\lambda_{N+1}(h),  \dots, \lambda_{N+K}(h))$, $h\in X_m$. They span a space $W\subset \R^K$ of dimension $\le m$.  We  assume that $\dim(W)=m$ (a slight notational adjustment has to be made if $\dim(W)<m$).   
 It follows from \eref{approx1} that each row ${\bf g}_i$ of $G$  can be approximated by a vector  from  $W$ in the $\ell_\infty$ norm to accuracy $d_m$, and therefore    in the $\ell_2$ norm to accuracy  $\sqrt{K}d_m$.
 Let $P$ be the orthogonal projection of $\R^K$ onto $W$.   Hence, we have  
 \be
 \label{also1}
\|{\bf g}_i-P{\bf g}_i\|_{\ell_2}\le \sqrt{K} d_m, \quad i=1,\dots , K.
 \ee
It also follows from \eref{entries} that 
\be
\nonumber
 \|P{\bf g}_i\|_{\ell_2}\le \|{\bf g}_i\|_{\ell_2}\le \left\{\sum_{j=1}^i\sigma^2_{N+j}\right\}^{1/2},
\ee
and therefore
\be
 \label{also2}
\sum_{i=1}^K \|P{\bf g}_i\|_{\ell_2}^2\le \sum_{i=1}^K\sum_{j=1}^i\sigma^2_{N+j}\leq K \sum_{i=1}^K \sigma^2_{N+i}.
\ee
Next, we apply Lemma \ref{L1} for this $G$ and $W$ and use estimates \eref{diag}, \eref{also1} and \eref{also2} to derive 
\begin{eqnarray}
\nonumber
 \gamma^{2K}\prod_{i=1}^K\sigma_{N+i}^2&\leq&\left\{\frac{K}{m}\sum_{i=1}^K \sigma^2_{N+i}\right\}^m\left\{ \frac{K^2}{K-m}d^2_m\right\}^{K-m}\\
\nonumber
&=& K^{K-m}\left(\frac{K}{m}\right)^m\left(\frac{K}{K-m}\right)^{K-m}
\left\{\sum_{i=1}^K\sigma_{N+i}^2\right\}^m d^{2(K-m)}_m\\
\nonumber
&\le &2^KK^{K-m}\left\{\sum_{i=1}^K\sigma_{N+i}^2\right\}^m d^{2(K-m)}_m ,
\end{eqnarray}
 and the proof is completed.
\hfill $\Box$

In analogy with Corollary \ref{C1}, we have the following special results for the weak
greedy algorithm in a general Banach space.
 \begin{cor}
 \label{C2}
 Suppose that $X$ is a Banach space.  For the weak greedy algorithm  with a constant $\gamma$, applied to a compact set $\cF$ contained in the unit ball of  $X$, the following  holds for $\sigma_n:= \sigma_n(\cF)_X$ and $d_n:=d_n(\cF)_X$, $n=1,2,\dots$,
 
 {\rm (i)} For any such compact set $\cF$ and $n\ge 1$, we have
 \be
 \label{C21}
 \sigma_n\le  \sqrt{2}\gamma^{-1}\min_{1\le m< n} n^{\frac{n-m}{2n}}\left\{\sum_{i=1}^n\sigma_i^2\right \}^{\frac{m}{2n}}d_m^{\frac{n-m}{n}}.
\ee 
In particular $\sigma_{2\ell}\le 2\gamma^{-1} \sqrt{\ell d_\ell}$, $\ell=1,2\dots$.

{\rm (ii)} If for $\alpha>0$, we have $d_n\le C_0n^{-\alpha}$, $n=1,2,\dots$, then for any $0<\beta<\min\{\alpha,1/2\}$, we have   $\sigma_n\le C_1n^{-\alpha+1/2+\beta}$, $n=1,2\dots$, with 

$$
C_1:=\max\left\{C_04^{4\alpha+1}\gamma^{-4}\left(\frac{2\beta+1}{2\beta}\right)^{\alpha},
\max_{n=1,\ldots,7}\{  n^{\alpha-\beta-1/2}\}\right\}.
$$

{\rm (iii)} If  for $\alpha>0$, we have $d_n\le C_0e^{-c_0n^{\alpha}}$, $n=1,2,\dots$, then    
$\sigma_n < \sqrt{2C_0}\gamma^{-1}\sqrt{n}e^{-c_1n^\alpha}$, $n=1,2\dots$, 
where $c_1=2^{-1-2\alpha}c_0$.  The factor $\sqrt{n}$ can be deleted by reducing  the constant  $c_1$.

 \end{cor}
 
 {\bf Proof:}  The proofs are  similar to those of  Corollary \ref{C1} except that we use \eref{T3} in place of \eref{T1}.
 
 (i) We take $N=0$, $K=n$, and any $1\le m< n$ in \eref{T3} and use the monotonicity of $(\sigma_n)_{ n\geq 0}$  to  obtain
 \be
 \label{C140}
 \sigma_n^{2n}\le 2^nn^{n-m} \gamma^{-2n}\left\{\sum_{i=1}^{ n} \sigma_i^2\right\}^m d_m^{2n-2m}. \ee
If we take a $2n$-th root of both sides, we arrive at \eref{C21}.  In particular, if $n=2\ell$ and $m=\ell$, we have
\be
\nonumber
\sigma_{2\ell}\le \sqrt{2}\gamma^{-1} (2\ell)^{1/4}\{\Sigma_{i=1}^{2\ell}\sigma_i^2\}^{1/4}\sqrt{d_\ell}\leq
\sqrt{2}\gamma^{-1} (2\ell)^{1/4}(2\ell)^{1/4}\sqrt{d_\ell}=2\gamma^{-1}\sqrt{\ell d_\ell},
\ee
where we have used the fact that all $\sigma_i\leq 1$.

 (ii)   It follows from the monotonicity of $(\sigma_n)_{n\geq 0}$ and \eref{T3} for $N=K=n$ and any $1\leq m<n$
 that
 \be
 \label{Cnew1}
  \sigma_{2n}
\le \sqrt{2n}\gamma^{-1}\sigma_{n}^{\delta}d_m^{(1-\delta)}, \quad \delta:=\frac{m}{n}.
 \ee
 Given our  $\beta$, we define  $m=:\lfloor\frac{2\beta}{2\beta+1}n\rfloor+1$ ($m<n$ for $n\geq2>2\beta+1$).   It follows that 
\be
\label{ds}
\delta=\frac{m}{n}\in \left(\frac{2\beta}{2\beta+1}, \frac{2\beta}{2\beta+1}+\frac{1}{n}\right).
\ee
We next prove (ii) by contradiction.  Suppose it is not true and $M$ is  the first value   
where $\sigma_M>C_1M^{ -\alpha+\beta+1/2}$.   Clearly, because of the definition of 
$C_1$,  and the fact that $\sigma_n\le 1$, we must have $M>7$. 
We  first  consider the case  $M=2n$, and therefore $n>3$.  From \eref{Cnew1} we have  
$$
C_1(2n)^{-\alpha+\beta+1/2} <\sigma_{2n}\leq \sqrt{2n}\gamma^{-1}C_1^\delta n^{\delta(-\alpha+\beta+1/2)}C_0^{1-\delta}(\delta n)^{-\alpha(1-\delta)},
$$
where we have used the fact that $\sigma_{n}\leq C_1n^{-\alpha+\beta+1/2}$ and $d_m\leq C_0m^{-\alpha}$. It follows that
$$
C_1^{1-\delta}<C_0^{1-\delta}2^{\alpha-\beta}\gamma^{-1}\delta^{-\alpha(1-\delta)}
n^{\frac{2\beta+1}{2}(\delta-\frac{2\beta}{2\beta+1})},
$$
and therefore
$$
C_1<C_02^{\frac{\alpha-\beta}{1-\delta}}\gamma^{-\frac{1}{1-\delta}}\delta^{-\alpha}
n^{\frac{2\beta+1}{2}\cdot\frac{\delta-\frac{2\beta}{2\beta+1}}{1-\delta}}.
$$
Since for $ n\geq 4>2(2\beta+1)$, we have
$$
\delta<\frac{4\beta+1}{2(2\beta+1)}<1,\quad \hbox {and therefore}\quad \frac{1}{1-\delta}<2(2\beta+1).
$$
 This gives
$$
{\frac{2\beta+1}{2}\cdot\frac{\delta-\frac{2\beta}{2\beta+1}}{1-\delta}}<\frac{(2\beta+1)^2}{n}, 
\quad \hbox{and thus}\quad n^{\frac{2\beta+1}{2}(\delta-\frac{2\beta}{2\beta+1})}<n^{\frac{(2\beta+1)^2}{n}}<2^{(2\beta+1)^2}.
$$
Then,  for $\beta<\min\{\alpha,1/2\}$ 
\begin{eqnarray}
\nonumber
C_1&<&C_02^{2(\alpha-\beta)(2\beta+1)}\gamma^{-2(2\beta+1)}\left(\frac{2\beta}{2\beta+1}\right)^{-\alpha}
2^{(2\beta+1)^2}\\
&<&C_02^{2(2\alpha+1)}\gamma^{-4}\left(\frac{2\beta+1}{2\beta}\right)^{\alpha}
<C_02^{2(4\alpha+1)}\gamma^{-4}\left(\frac{2\beta+1}{2\beta}\right)^{\alpha},
\end{eqnarray}
 which is the desired contradiction. Likewise, if $M=2n+1$ (since $M>7$, we have $n>3$), 
 for $-\alpha+\beta+1/2<0$ (which is the meaningful case),
\begin{eqnarray}
\nonumber
C_12^{-\alpha+\beta+1/2}(2n)^{-\alpha+\beta+1/2}&<&C_1(2n+1)^{-\alpha+\beta+1/2}<\sigma_{2n+1}\leq \sigma_{2n}\\
&<& 
\nonumber
\sqrt{2n}\gamma^{-1}C_1^\delta n^{\delta(-\alpha+\beta+1/2)}C_0^{1-\delta}(\delta n)^{-\alpha(1-\delta)},
\end{eqnarray}
and following the same argument as above we get
$$
C_1<C_02^{-\frac{1}{2(1-\delta)}}2^{2\frac{\alpha-\beta}{1-\delta}}
\gamma^{-2(2\beta+1)}\left(\frac{2\beta}{2\beta+1}\right)^{-\alpha}
2^{(2\beta+1)^2}<C_02^{2(4\alpha+1)}\gamma^{-4}\left(\frac{2\beta+1}{2\beta}\right)^{\alpha},
$$
where we have used that $2^{-\frac{1}{2(1-\delta)}}<1$, and the proof  is completed.

 (iii) From (i), we have
 \be
\nonumber
 \sigma_{2n+1} \le \sigma_{2n} \le 
2\gamma^{-1}\sqrt{nd_n}\le 2\gamma^{-1}\sqrt{C_0}\sqrt{n}e^{-\frac{c_0}{2}n^\alpha}
< \sqrt{2C_0}\gamma^{-1}\sqrt{2n+1}e^{-c_02^{-1-2\alpha}(2n+1)^\alpha}, 
 \ee
 from which (iii) easily follows.  \hfill $\Box$
 
 \section{Lower bounds in a Banach space}  
 
 It is natural to ask whether the factor $\sqrt{n}$ is necessary when proving results in a Banach space.  
 Here, we shall provide examples which show that a  loss of this type  is indeed necessary.   However,
 as   it will be seen,  there is still a small gap between what we have proved for direct estimates and what the
 examples  below provide.

Let us begin by considering the space $X:=\ell_\infty (\N\cup\{0\})$ equipped with its usual norm.  We consider a monotone decreasing sequence $x_0\ge x_1\ge x_2\ge \cdots$
of positive real numbers   which converge to zero and define  
\be
\nonumber
f_j:=x_je_j,\quad j=0,1,\dots,
\ee
where $e_j$, $j=0,1,\dots$ are the usual coordinate vectors in $\R^{\N\cup\{0\}}$.    Let $\cF:=\{f_0,f_1,\dots \}$.  From the monotonicity of the $x_j$'s, the greedy algorithm for $\cF$ in $X$ can  choose the elements from $\cF$ in order $f_0,f_1,\dots$.  Hence, 
\be
\nonumber
\sigma_j=\sigma_j(\cF)_{X}=x_j,\quad j\ge 0.
\ee
We   want to give an upper bound for the Kolmogorov width of $\cF$.  For this, we shall use the following
result (see (7.2) of Chapter 14  in  \cite{LGM}) on $s$-widths  of the unit ball $b_1^m$ of $\ell_1^m$ in $\ell_\infty^m$:
\be
\label{pajor}
d_s(b_1^m)_{X}\le C \left\{\ln (m/s)\right\}^{1/2} s^{-1/2}, \quad 1\le s\le m/2.
\ee
Let us now  define  the  sequence  $\{x_j\}_{j\geq 0}$ so that in position $2^{k-1}\le j\le 2^{k}-1$ it has the constant value $2^{-k\alpha}$, for $k=0,1,\dots$,  where $\alpha>1/2$.  It follows that,
\be
\nonumber
\sigma_n(\cF)_{X} =\cO( n^{-\alpha}), \quad n=1,2,\dots.
\ee

 We shall now bound the $N$-width of $\cF$ when $N=2^{n+1}$ by constructing a good space $X_N$ 
  of dimension $\leq N$ for approximating $\cF$.   The space $X_N$ will be the span of a set $E$  of at most $N$ vectors.  First, we  place in $E$ all of the vectors, $e_1,\dots,e_{2^n}$.  Next, for each $k=1,\dots n$, we  use \eref{pajor}  to choose a    basis for the  space of dimension
$2^{n-k} $ whose vectors are supported on $[2^{n+k},2^{n+k+1}-1]$ and this space approximates each of the $f_j$, $j=2^{n+k},\dots,2^{n+k+1}-1$, in $X$  to accuracy $C_02^{-(n+k)\alpha} \sqrt{k}2^{-(n-k)/2}$.
 We place these basis vectors in $E$.   Notice that $|x_j|\le 2^{-2n\alpha}$ for $j\ge 2^{2n}$. 
 This means that for the space $X_N:=\span(E)$  with  dimension $\le N$ we have
\begin{eqnarray}
\nonumber
d_N(\cF)_{X} &\le& \dist (\cF,X_N)_X\le \max  \left\{2^{-2n\alpha}, \max_{1\le k\leq n}C_02^{-(n+k)\alpha} 2^{-(n-k)/2} \sqrt{
k}\right\}\nonumber \\
&=& \max  \left\{2^{-2n\alpha}, C_02^{-n(\alpha+1/2)}\cdot\max_{1\le k\leq n}2^{-k(\alpha-1/2)}\sqrt{ k}\right\}
\le C_1 2^{-n(\alpha+1/2)}, \quad  \alpha>1/2.
\nonumber
\end{eqnarray}
From the monotonicity of $(d_n(\cF)_X)_ {n\geq 0}$, we obtain that
\be 
\nonumber
d_n(\cF)_{X}  \le C_2 n^{-\alpha-1/2},\quad n=1,2,\dots.
\ee

This example shows that the factor $\sqrt{n}$ which appears in (ii) of Corollary \ref{C2} can in general not be removed.


 \vskip .1in
\noindent
Ronald DeVore, Department of Mathematics, Texas A\&M University,
College Station, TX,
  rdevore@math.tamu.edu
\vskip .1in
\noindent
Guergana Petrova, Department of Mathematics, Texas A\&M University,
College Station, TX,    gpetrova@math.tamu.edu
\vskip .1in
\noindent
Przemyslaw Wojtaszczyk,  Institute  of Applied Mathematics, and Interdisciplinary Centre for Mathematical and Computational Modelling, University of Warsaw,  Warsaw, Poland, wojtaszczyk@mimuw.edu.pl

\end{document}